\documentclass[10pt,twoside]{article}
\usepackage{amsmath}
\usepackage{amssymb}
\usepackage{amsthm}
\usepackage{Latex-document}
\def\XXint#1#2#3{{\setbox0=\hbox{$#1{#2#3}{\int}$}
     \vcenter{\hbox{$#2#3$}}\kern-.5\wd0}}

\newcommand{\bd}[1]{\partial #1}
\newcommand{\cl}[1]{\overline{#1}}
\newcommand{\pd}[2]{\frac{\partial #1}{\partial #2}}
\newcommand{\R}{\mathbb R}
\DeclareMathOperator{\BMO}{BMO} \DeclareMathOperator{\supp}{supp}
\DeclareMathOperator{\VMO}{VMO} \theoremstyle{plain}
\newtheorem{thm}{Theorem}
\newtheorem*{enthm}{Theorem \thmn}

\theoremstyle{remark}

\numberwithin{equation}{thm}

\title{\bf Harmonic Measure and
 ``Locally Flat''\vskip -2mm Domains\thanks{Partially
supported by the NSF.}  \vskip 6mm}

\author{Carlos E.\ Kenig\vspace*{-0.5cm}\thanks{University of Chicago, Department of Mathematics, Chicago,
IL 60637, USA.  E-mail:  cek@math.uchicago.edu  }}

\date{\vspace{-8mm}}

\begin{document}

\maketitle

\thispagestyle{first} \setcounter{page}{701}

\begin{abstract}\vskip 3mm
We will review work with Tatiana Toro yielding a characterization
of those domains for which the harmonic measure has a density
whose logarithm has vanishing mean oscillation. \vskip 4.5mm


\end{abstract}

\vskip 12mm

In this lecture, I will describe a series of joint works
 with Tatiana Toro on the relationship between regularity
properties of harmonic measure and Poisson kernels, and regularity
properties of the underlying domains.  Thus, consider a domain
$\Omega \subseteq \R^{n + 1}$ and the solution to the classical
Dirichlet problem:
\begin{equation}
\tag{DP}
\begin{cases}
\Delta u = 0 & \text{in }\Omega                 \\
u\bigr|_{\bd{\Omega}} = f \in C_b(\bd{\Omega}),
\end{cases}
\end{equation}
$u \in C_b(\cl{\Omega})$, where $C_b$ is the class of bounded
continuous functions.  The maximum principle and the Riesz
representation theorem yield the formula $$u(X_*) =
\int_{\bd{\Omega}} f(Q)d\omega^{X_*}(Q),\quad X_* \in \Omega,$$
and the family of positive Borel probability measures
$\{d\omega^{X_*}\}$ is called harmonic measure.  We sometimes fix
$X_* \in \Omega$ and write $d\omega = d\omega^{X_*}$.  Note that,
if $\Omega$ is a smooth domain, then $d\omega^{X_*}(Q) =
\pd{G}{\vec{n}_Q}(Q, X_*)d\sigma(Q)$, where $G$ is the Green's
function for $\Omega$, $d\sigma$ is surface measure, and
$\pd{}{\vec{n}_Q}$ denotes differentiation along the outward unit
normal.  When $\Omega$ is unbounded and $v$ is a minimal harmonic
function in $\Omega$ with $v\bigr|_{\bd{\Omega}} \equiv 0$, we
define $d\omega^\infty$, harmonic measure with pole at infinity,
to be the measure satisfying $$\int_{\bd{\Omega}} \varphi
d\omega^\infty = \int_\Omega v\Delta \varphi,\quad\text{for
}\varphi \in C_0^\infty(\Omega).$$  The existence and uniqueness
of $v$ and $\omega^\infty$ (modulo multiplicative constants) can
be established, for instance, when $\Omega$ is an unbounded NTA
(non-tangentially accessible) domain
%%In the original, this is just `NTA'.  Since the
abbreviation is clarified later, it might as well be here. (see
\cite{KT1} for details).  For example, if $\Omega = \R^{n + 1}_+ =
\{(x, t) : t > 0\}$, then $v(x, t) = t$ and $d\omega^\infty = dx$
on $\R^n$.  The work I will describe originated from trying to
understand, as $\alpha \to 0$, the classical theorem of Kellogg,
which shows that, if $\Omega$ is of class $C^{1, \alpha}$, $0 <
\alpha < 1$, then $d\omega = k\,d\sigma$ with $\log k \in
C^{\alpha}$; and its ``converse'', the free boundary regularity of
Alt-Caffarelli \cite{AC}, which states that, if $\Omega$ satisfies
certain necessary weak conditions (to be more fully explained
later) and $d\omega = k\,d\sigma$ with $\log k \in C^\alpha$, then
$\Omega$ must be of class $C^{1, \alpha}$.

To motivate our results, we recall real variable characterizations
of $C^{1, \alpha}$ and $C^\alpha$:
\begin{multline}
\tag*{(I)$_\alpha$}
\varphi \in C^{1, \alpha}(\R^n)\text{ (}0 < \alpha < 1\text{)} \Leftrightarrow \forall r > 0, x_0 \in \R^n\text{, there exists an affine function} \\
L_{r, x_0}\text{ on }\R^n\text{ such that }\frac{|\varphi(x) -
L_{r, x_0}(x)|}{r} \leq C r^\alpha\text{ for }|x - x_0| < r.
\end{multline}
When $\alpha = 0$, this condition is equivalent to the Zygmund
class condition $\varphi \in \Lambda_*$, i.e.,
$$\frac{|\varphi(x + h) + \varphi(x - h) -
2\varphi(x)|}{|h|} \leq C.$$  For us, when $\alpha = 0$, the
$\lambda_*$ class will also be relevant, where $\varphi \in
\lambda_*$ if $\varphi \in \Lambda_*$ and, in addition, the ratio
described above tends to $0$ as $\alpha \to 0$.
\begin{equation}
\tag*{(II)$_\alpha$} h \in C^\alpha \Leftrightarrow \sup_{r > 0}
\frac{1}{r^\alpha}\text{av}_{B_R} |h - h_{B_r}| \leq C,
\end{equation}
where $\text{av}_A$ denotes the average over the set $A$ and $B_r$
any ball of radius $r$.  When $\alpha = 0$, this becomes the BMO
space of John-Nirenberg \cite{JN}, but we will be more interested
in $\VMO$, where $h \in \VMO$ if $h \in \BMO$ and in addition
$\text{av}_{B_r} |h - h_{B_r}| \xrightarrow{r \to 0} 0$. Note that
$\VMO$ plays the role {\em vis-\`a-vis} $\BMO$ that continuous
functions play {\em vis-\`a-vis} $L^\infty$.

We start out by giving our geometric analogue of (I)$_0$: We say
that $\Omega \subseteq \R^{n + 1}$ is $\delta$-Reifenberg flat if
it has the separation property (a quantitative connectivity
property) (see \cite{KT1} for details), and, for all compact $K
\subseteq\subseteq \R^{n + 1}$, there exists $R_K > 0,$ such that,
for $0 < r < R_k$ and $Q \in \bd{\Omega} \cap K$, there exists an
$n$-dimensional plane $L(r, Q)$ passing through $Q$ such that
$$\frac{1}{r}D\bigl[B(r, Q) \cap \bd{\Omega}, B(r, Q) \cap
L(r, Q)\bigr] \leq \delta,$$ where $D$ denotes Hausdorff distance.
Note that this is a significant condition only for $\delta < 1$.
We will always assume $\delta < \frac{1}{4\sqrt{2}}$.  We say that
$\Omega$ is Reifenberg vanishing if, as $r \to 0$, we can take
$\delta \to 0$.  For instance, the domain above the graph of a
$\lambda_*$ function is Reifenberg vanishing.  In general,
Reifenberg vanishing domains are not local graphs; they do not
have tangent planes or a ``surface measure''.  This class of
domains was introduced by Reifenberg \cite{R} in his study of the
Plateau problem for minimal surfaces in higher dimensions.

In order to state our analogue of Kellogg's theorem in this setting, we need to introduce ``multiplicative'' analogues %%The spelling switched here.  Earlier it was `analog'.
of (I)$_0$.  A measure $\mu$, supported on $\bd{\Omega}$, is doubling if, $\forall K \subseteq\subseteq \R^{n + 1}$, there exists $R_K > 0$ such that, if $0 < r < R_K$, then $$\mu\bigl(B(2r, Q) %%This is somewhat erratic.  Sometimes it's $B(2r, Q)$ and sometimes $B(Q, 2r)$.
\cap \bd{\Omega}\bigr) \leq C\,\mu\bigl(B(r, Q) \cap
\bd{\Omega}\bigr).$$  Such a $\mu$ is called asymptotically
optimal doubling (see \cite{}, \cite{}) for details) if it is
doubling and $$\lim_{r \to 0} \inf_{Q \in \bd{\Omega} \cap K}
\frac{\mu(B(\tau r, Q) \cap \bd{\Omega})}{\mu(B(r, Q) \cap
\bd{\Omega})} = \lim_{r \to 0} \sup_{Q \in \bd{\Omega} \cap K}
\frac{\mu(B(\tau r, Q) \cap \bd{\Omega})}{\mu(B(r, Q) \cap
\bd{\Omega})} = \tau^n,$$ for $0 < \tau < 1$, $K
\subseteq\subseteq \R^n$.  For example, if $\Omega$ is of class
$C^{1, \alpha}$ and $d\sigma$ denotes surface measure, then
$\sigma\bigl(B(r, Q) \cap \bd{\Omega}) = \alpha_nr^n + O(r^{n +
\alpha})$, $Q \in \bd{\Omega}$, and hence $\sigma$ is
asymptotically optimal doubling.  If $\log k \in C^\alpha$, then
the same is true for $d\omega = k\,d\sigma$.  Our analog of
Kellog's theorem is:
\begin{thm}
\label{Thm1} {\rm (\cite{KT3})}  If $\Omega$ is a Reifenberg
vanishing domain, then $\omega$ ($\omega^\infty$) is
asymptotically optimal doubling.
\end{thm}

The proof uses the fact that $\delta$-Reifenberg flat domains are
NTA domains (\cite{JK1}, \cite{KT3}).  One then uses the theory of
the boundary behavior of harmonic functions on NTA domains
(\cite{JK1}) and comparisons to half-planes, using the Reifenberg
vanishing condition and the maximum principle.

To understand a possible converse to Theorem \ref{Thm1},
 we recall a geometric measure theory (GMT) problem, first
posed by Besicovitch:  let $\mu$ be a positive Radon measure on
$\R^{n + 1}$ such that, for each $Q \in \Sigma$ ($\Sigma$ the
support of $\mu$) and each $r > 0$, we have
\begin{equation}
\tag{B} \mu\bigl(B(r, Q)\bigr) = \alpha r^n,\quad \alpha > 0\text{
fixed.}
\end{equation}
Then, what can be said about $\mu$?  Clearly, if $d\mu = dx$ on
$\R^n \subseteq \R^{n + 1}$, then (B) holds. Nevertheless, in
1987, D.\ Preiss found the following interesting example:  let
$\Sigma_C$ be the light cone $x_4^2 = x_1^2 + x_2^2 + x_3^2$, and
$d\mu = d\sigma_{\Sigma_C}$ its surface measure.  Then $\mu$
satisfies (B).  Moreover, the general case of (B) is settled by
the following remarkable theorem of Kowalski-Preiss \cite{KoP}.
\vskip 3mm

\noindent{\bf Theorem.}  {\rm (\cite{KoP})}  {\it Let $\mu$ be a
non-zero measure with property \textup{(B)}, and put $\Sigma =
\supp \mu \subseteq \R^{n + 1}$. If $n = 1, 2$, then $\Sigma =
\R^n$.  If $n \geq 3$, then either $\Sigma = \R^n$ or $\Sigma =
\Sigma_C \otimes \R^{n - 3}$, modulo rigid motions.} \vskip 3mm

The connection of the Preiss example to our problem comes from the
fact (\cite{KT1}) that, if $\Omega = \{x_4^2 < x_1^2 + x_2^2 +
x_3^2\}$, $\Omega \subseteq \R^4$, then $d\omega^\infty =
d\sigma_{\Sigma_C}$ (separation of variables) and, by Preiss's
result, $\omega^\infty$ is asymptotically optimal doubling, but,
of course, $\Omega$ is not Reifenberg vanishing, since it is
$\frac{1}{4\sqrt{2}}$-Reifenberg flat, and no better.  Our
converse to Theorem \ref{Thm1} is now:
\begin{thm}
\label{Thm2} {\rm (\cite{KT1})}  Assume that $\Omega \subseteq
\R^{n + 1}$ is an NTA, and that $\omega$ ($\omega^\infty$) is
asymptotically optimal doubling.  If $n = 1, 2$, then $\Omega$ is
Reifenberg vanishing.  If $n \geq 3$ and $\Omega$ is
$\delta$-Reifenberg flat, $\delta < \frac{1}{4\sqrt{2}}$, then
$\Omega$ is Reifenberg vanishing.
\end{thm}

This is in fact a GMT result.  It remains valid if $\omega$
($\omega^\infty$) is replaced by any asymptotically optimal
doubling measure $\mu$ with support $\bd{\Omega}$.  The idea of
the proof is to use a ``blow-up'' argument to reduce matters to
the Kowalski-Preiss theorem.  Further GMT results along these
lines, also in the higher codimension case, were obtained by
David-Kenig-Toro \cite{DKT}.

We now turn to the results motivated by (II)$_0$.  Note that the
unit normal $\vec{n}$ satisfies $|\vec{n}| = 1$, and so the BMO
condition on it is automatic, but the VMO condition is not.  To
put our work in perspective, we recall some of the history of the
subject.

A domain $\Omega \subseteq \R^{1 + 1} = \R^2$ is called a
chord-arc domain is $\bd{\Omega}$ is locally rectifiable, and,
whenever $Q_1, Q_2 \in \bd{\Omega}$, we have $\ell\bigl(s(Q_1,
Q_2)\bigr) \leq C|Q_1 - Q_2|$, where $\ell$ denotes length and
$s(Q_1, Q_2)$ is the shortest arc between $Q_1$ and $Q_2$.
$\Omega$ is called vanishing chord-arc if, in addition, as $Q_1
\to Q_2$, the ratio $\frac{\ell(s(Q_1, Q_2))}{|Q_1 - Q_2|}$ tends
to $1$, uniformly on compact sets.  The first person to study
harmonic measure on chord-arc domains in the plane was Lavrentiev
(\cite{L}), who proved:\vskip 3mm

\noindent{\bf Theorem.}  {\rm (\cite{L})} {\it If $\Omega
\subseteq \R^{1 + 1}$ is chord-arc, then $d\omega = k\,d\sigma$
with $\log k \in \BMO(d\sigma)$.  (In fact, $\omega \in
A_\infty(d\sigma)$, the Muckenhoupt class \cite{GCRdeF}.)}\vskip
3mm

For vanishing chord-arc domains in the plane, Pommerenke \cite{P}
proved:\vskip 3mm

\noindent{\bf Theorem.}  {\rm (\cite{P})}  {\it Suppose that
$\Omega$ is a chord-arc domain in $\R^{1 + 1}$.  Then $\Omega$ is
vanishing chord-arc if and only if $d\omega = k\,d\sigma$ with
$\log k \in \VMO(d\sigma)$.}\vskip 3mm

These results were obtained using function theory, so their proofs don't generalize to higher dimensions.  In higher dimensions, the first breakthrough came in the celebrated theorem of B. Dahlberg \cite{D}, who showed that, if $\Omega \subseteq \R^{n + 1}$ is a Lipschitz domain, then $d\omega = k\,d\sigma%%In the text this also seems to be $d\omega$.
$ with $\log k \in \BMO$ (in fact, $\omega \in
A_\infty(d\sigma)$).  One direction of Pommerenke's result was
extended to higher dimensions by Jerison-Kenig \cite{JK2}, who
showed that, if $\Omega$ is a $C^1$ domain, then $\log k \in
\VMO$.  (In general, note that $\Omega$ is of class $C^1$ need not
imply that $\log k$ is continuous.)  In order to explain our
results and to clarify the connection with condition (II)$_0$, we
need to introduce some terminology. A domain $\Omega \subseteq
\R^{n + 1}$ will be calld a chord-arc domain if it is an NTA
domain (see \cite{JK1}) of locally finite perimeter (see
\cite{EG}) and its boundary is Ahlfors regular, i.e., the surface
measure $\sigma$ (which is Radon measure on $\bd{\Omega}$, by the
assumption of locally finite perimeter) satisfies the inequalities
$$C^{-1}r^n \leq \sigma\bigl(B(r, Q) \cap \bd{\Omega}\bigr)
\leq C r^n$$ (for $Q \in K \cap \bd{\Omega}$, $K
\subseteq\subseteq \R^{n + 1}$ and small $r$; or, if $\Omega$ is
an unbounded NTA, for all $Q \in \bd{\Omega}$ and $r > 0$).  A
fundamental result of David-Jerison \cite{DJ} and Semmes
\cite{Se2} is that Dahlberg's theorem extends to this case, i.e.,
that $d\omega = k\,d\sigma$ with $\log k \in \BMO$, and, in fact,
$\omega \in A_\infty(d\sigma)$.

We say that $\Omega \subseteq \R^{n + 1}$ is a
``$\delta$-chord-arc domain'' if $\Omega$ if $\delta$-Reifenberg
flat, $\Omega$ is of locally finite perimeter, the boundary of
$\Omega$ is Ahlfors regular and the BMO norm of the unit normal
$\vec{n}$ is bounded by $\delta$.  We say that $\Omega$ is
``vanishing chord-arc'' if, in addition, it is Reifenberg
vanishing and $\vec{n} \in \VMO(d\sigma)$.\vskip 3mm

\noindent{\bf Remark. }
S.\ Semmes \cite{Se1} has proved that, if $\Omega$ is a $\delta$-chord-arc domain (under the definition used above), then $$(1 + \sqrt{\delta})^{-1}\alpha_nr^n \leq \sigma\bigl(B(r, Q) \cap \bd{\Omega}\bigr) \leq (1 + \sqrt{\delta})\alpha_nr^n,$$ where $\alpha_n$ is the volume of the unit ball in $\R^n$ and $\delta \leq \delta_n$%%What's $\delta_n$?  It's mentioned several times, but never defined.
.  Moreover, a combination of the results in \cite{Se1} and
\cite{KT1} shows that, if $\Omega$ is a $\delta$-Reifenberg flat
domain which is of locally finite perimeter, and for which
$\sigma\bigl(B(r, Q) \cap \bd{\Omega}) \leq (1 +
\delta)\alpha_nr^n$, then the BMO norm of $\vec{n}$ is bounded by
$C\sqrt{\delta}$ for $\delta < \delta_n$.  Thus, the two notions
introduced of ``vanishing chord-arc'' domains in the plane are the
same, and a domain is vanishing chord-arc exactly when it is of
locally finite perimeter, has an Ahlfors regular boundary, it is
Reifenberg vanishing and satisfies $\vec{n} \in \VMO$. \vskip 3mm

Our potential-theoretic result, which extends the work of
Jerison-Kenig \cite{JK2}, is
\begin{thm}
\label{Thm2'}%%In the text, this was another Theorem 2.
{\rm (\cite{KT3})}  If $\Omega$ is a vanishing chord-arc domain,
then $\omega$ ($\omega^\infty$) has the property that $d\omega =
k\,d\sigma$ ($d\omega^\infty = h\,d\sigma$) with $\log k \in \VMO$
($\log h \in \VMO$).
\end{thm}

This was proved by a combination of real variable arguments,
potential-theoretic arguments, and the estimates in \cite{JK2}.

In order to understand possible converses of this, extending the
work of Pommerenke to higher dimensions, we will recall precisely
the Alt-Caffarelli \cite{AC} result which we alluded to earlier.
In the language that we have introduced, their local regularity
theorem can be stated as follows:\vskip 3mm

\noindent{\bf Theorem.}  {\rm (\cite{AC})} {\it Let $\Omega$ be a
set of locally finite perimeter whose boundary is Ahlfors regular.
Assume that $\Omega$ is $\delta$-Reifenberg flat, $\delta <
\delta_n$. Suppose that $d\omega = k\,d\sigma$ with $\log k \in
C^\alpha(\bd{\Omega})$ ($0 < \alpha < 1$).  Then $\Omega$ is a
$C^{1, \alpha}$ domain.}\vskip 3mm

The reason for this being a free boundary regularity result is
that, in the case when $\Omega$ is unbounded and $d\omega =
d\omega^{\infty}$, $d\omega^\infty = h\,d\sigma$, then $v > 0$ in
$\Omega$, $v\bigr|_{\bd{\Omega}} \equiv 0$, $\Delta v = 0$ in
$\Omega$ and $h = \pd{v}{\vec{n}}$.  Thus, knowledge of the
regularity of the Cauchy data of $v$ ($v\bigr|_{\bd{\Omega}}$,
$\pd{v}{\vec{n}}\bigr|_{\bd{\Omega}}$) yields regularity of
$\bd{\Omega}$ (or of $\vec{n}$, the normal).

The first connection between the above Theorem and the work of
Pommerenke was made by Jerison \cite{J}, who was also the first to
formulate the higher-dimensional analogues of Pommerenke's theorem
as end-point estimates as $\alpha \to 0$ in the Alt-Caffarelli
theorem. Jerison's theorem in \cite{J} states that, if $\Omega$ is
a Lipschitz domain and $d\omega = k\,d\sigma$ with $\log k$
continuous, then $\vec{n} \in \VMO$.  There is an error in Lemma 4
of Jerison's paper.  Nonetheless, in \cite{KT1} we made strong use
of the ideas in \cite{J}.  In the more recent version of our
results \cite{KT2}, we bypass this approach.

Before stating our result, it is useful to classify the
assumptions in the Alt-Caffarelli theorem.  For this, we recall
some examples: \vskip 2mm

\noindent {\bf Examples.}  When $n = 1$, Keldysh-Lavrentiev
\cite{KL} (see also \cite{Du}) constructed domains in $\R^{1 +
1}$ with locally rectifiable boundaries which (\cite{Du}) can be
taken to be Reifenberg vanishing and for which $d\omega =
d\sigma$, i.e., $k \equiv 1$, but which are not very smooth.  For
instance, they fail to be chord-arc.  These domains do not, of
course, have Ahlfors regular boundaries.  When $n = 2$,
Alt-Caffarelli constructed a double cone $\Gamma$ in $\R^3$ such
that, for $\Omega$ the domain outside the cone, $d\omega^\infty =
d\sigma$, i.e., $k \equiv 1$.  This is of course not smooth near
the origin, the problem being that, while $\Omega$ is NTA and
$\bd{\Omega}$ is Ahlfors regular, $\Omega$ is not
$\delta$-Reifenberg flat for small $\delta$.  When $n = 3$, the
Preiss cone we saw before exhibits the same behavior.

Our first result was:
\begin{thm}
\label{Thm3} {\rm (\cite{KT1})}  Assume that $\Omega \subseteq
\R^{n + 1}$ is $\delta$-chord-arc, $\delta \leq \delta_n$, that
$\omega$ ($\omega^\infty$) is asymptotically optimal doubling and
that $\log k \in \VMO$ ($\log h \in \VMO$).  Then $\vec{n} \in
\VMO$ and $\Omega$ is vanishing chord-arc.
\end{thm}

Notice, however, that, when comparing the hypothesis of Theorem
\ref{Thm3} to the Alt-Caffarelli theorem two things are apparent :
first, we are making the additional assumption that $\omega$ is
asymptotically optimal doubling, and hence, in light of Theorem
\ref{Thm2}, $\Omega$ is Reifenberg vanishing.  Next, the
``flatness'' assumption in the Alt-Caffarelli theorem is
$\delta$-Reifenberg flatness, while in Theorem \ref{Thm3} we make
the {\em a priori} assumption that, in addition, the BMO norm of
$\vec{n}$ is
smaller than $\delta$.  R%%In the text, this is ``Finally, r
...''.  This does not make much sense. ecently we have developed a
new approach which has removed these objections.  We have:
\begin{thm}
\label{Thm4} {\rm (\cite{KT2})}  Let $\Omega$ be a set of locally
finite perimeter whose boundary is Ahlfors regular.  Assume that
$\Omega$ is $\delta$-Reifenberg flat, $\delta < \delta_n$.
Suppose that $d\omega = k\,d\sigma$ ($d\omega^\infty =
h\,d\sigma$) with $\log k \in \VMO(d\sigma)$ ($\log h \in
\VMO(d\sigma)$).  Then $\vec{n} \in \VMO(d\sigma)$ and $\Omega$
is a vanishing chord-arc domain.
\end{thm}

Note that Theorems \ref{Thm2'} and \ref{Thm4} together give a
complete characterization of the vanishing chord-arc domains in
terms of their harmonic measure, in analogy with Pommerenke's
$2$-dimensional result, thus answering a question posed by Semmes
\cite{Se2}.

Our technique for the proof of Theorem \ref{Thm4} %%The text says Theorem 5 (really 6), but that hasn't occurred yet (and keeping the reference would make things circular).
is to use a suitable ``blow-up'' to reduce matters to the
following version of the ``Liouville theorem'' of Alt-Caffarelli
(\cite{AC}, \cite{KT5}):
\begin{thm}
\label{Thm5} {\rm (\cite{AC}, \cite{KT5})}  Let $\Omega$ be a set
of locally finite perimeter whose boundary is (unboundedly)
Ahlfors regular.  Assume that $\Omega$ is an unbounded
$\delta$-Reifenberg flat domain, $\delta < \delta_n$.  Suppose
that $u$ and $h$ satisfy:  $$\begin{cases}
\Delta u = 0                   & \text{in }\Omega \\
u > 0                          & \text{in }\Omega\
u\bigr|_{\bd{\Omega}} \equiv 0
\end{cases}$$ and $$\int_\Omega u\Delta\varphi = \int_{\bd{\Omega}} \varphi h\,d\sigma,\quad\text{for }\varphi \in C_0^\infty(\R^{n + 1}).$$  Suppose that $\sup_{x \in \Omega} |\nabla u(x)| \leq 1$ and $h(Q) \geq 1$ for ($d\sigma$-)a.e. $Q$ on $d\Omega$.  Then $\Omega$ is a half-space and $u(x, x_{n + 1}) = x_{n + 1}$.
\end{thm}

This allows us to prove the crucial blow-up result, which we now
describe.  Let $\Omega$ be as in Theorem \ref{Thm4}, and assume in
addition that $\Omega$ is unbounded.  Suppose $d\omega^\infty =
h\,d\sigma$ with $\log h \in \VMO(d\sigma)$, and let $u$ be the
associated harmonic function.  Let $Q_i \in \bd{\Omega}$ and
assume that $Q_i \to Q_\infty \in \bd{\Omega}$ as $i \to \infty$
(without loss of generality, $Q_\infty = 0$).  Let $\{r_i\}_{i =
1}^\infty$ be a sequence of positive numbers tending to $0$, and
put
\begin{gather*}
\Omega_i = \frac{1}{r_i}(\Omega - Q_i), \qquad \bd{\Omega}_i = \frac{1}{r_i}(\bd{\Omega} - Q_i), \\
u_i(X) = \frac{1}{r_i\text{av}_{B(r_i, Q_i)} h\,d\sigma}u(r_iX +
Q_i)\text{\ \ \ and\ \ \ }d\omega_i^\infty = h_i(Q)d\sigma_i(Q),
\end{gather*}
where $h_i(Q) = \frac{1}{\text{av}_{B(r_i, Q_i)} h\,d\sigma}h(r_iQ
+ Q_i)$.  Then:
\begin{thm}
\label{Thm6} There exists a subsequence of $\{\Omega_i\}$ (which
we will call again $\{\Omega_i\}$) satisfying:
\begin{gather}
\label{6.1}  \Omega_i \to \Omega_\infty\text{ in the Hausdorff distance sense, uniformly on compact sets;} \\
\label{6.2}  \bd{\Omega}_i \to \bd{\Omega_\infty}\text{ in the Hausdorff distance sense, uniformly on compact sets;} \\
\label{6.3}  u_i \to u_\infty\text{ uniformly on compact sets} \\
\intertext{and} \label{6.4}  \begin{cases}
\Delta u_\infty = 0 & \text{in }\Omega_\infty      \\
u_\infty = 0        & \text{in }\bd{\Omega_\infty} \\
u_\infty > 0        & \text{in }\Omega_\infty.
\end{cases}
\end{gather}
Furthermore
\begin{gather}
\label{6.5}  \omega_i \rightharpoonup \omega_\infty  \\
\intertext{and} \label{6.6}  \sigma_i \rightharpoonup
\sigma_\infty,
\end{gather}
weakly as Radon measures.  Here, $\sigma_\infty = {\mathcal H}^n \lfloor %%This isn't the right symbol, but that symbol doesn't seem to be available.
\bd{\Omega_\infty}$ and $\omega_\infty$ denotes the harmonic
measure of $\Omega_\infty$ with pole at $\infty$ (corresponding to
$u_\infty$).  Moreover,
\begin{gather}
\label{6.7}  \sup_{Z \in \Omega_\infty} |\nabla u_\infty(Z)| \leq
1 \\ \intertext{and} \label{6.8} h_\infty(Q) =
\frac{d\omega_\infty}{d\sigma_\infty}(Q) \geq 1\quad\text{ for
}{\mathcal H}^n\text{-a.e. }Q \in \bd{\Omega_\infty}.
\end{gather}
\end{thm}
Since $\log h \in \VMO(\bd{\Omega})$, the average $\text{av}_{B(r,
Q)} h\,d\sigma$ is close to the value of $\log h$ in a
proportionally large subset of $B(r, Q) \cap \bd{\Omega}$.  This
remark allows us to conclude that \eqref{6.6} holds, which is
crucial to the application and which fails in general under just
\eqref{6.1} and \eqref{6.2}.

As an immediate application of Theorems \ref{Thm5} and \ref{Thm6},
we obtain that $\Omega_\infty$ is a half-plane. This already
establishes that $\Omega$ is Reifenberg vanishing in Theorem
\ref{Thm4}.  To establish that $\vec{n}$ is in $\VMO$, we assume
otherwise, and obtain $Q_i \to Q_\infty$, $r_i \to 0$, such that
$\text{av}_{B(r_i, Q_i)} |\vec{n} - \vec{n}_{B(r_i,
Q_i)}|^2d\sigma \geq \ell^2$, $\ell > 0$.  We consider the
corresponding blow-up sequence, and let $\vec{e}_{n + 1}$ be the
direction perpendicular to $\bd{\Omega_\infty}$.  By the
divergence theorem and \eqref{6.1} and \eqref{6.2}, we have for
$\varphi \in C_0^\infty(\R^{n + 1})$ that $$\lim_{i \to \infty}
\int_{\bd{\Omega}_i} \varphi\langle\vec{n}_i, \vec{e}_{n +
1}\rangle d\sigma_i = \int_{\R^n \times \{0\}} \varphi\,dx$$ and
hence $$\lim_{i \to \infty} \Bigl\{\int_{\bd{\Omega}_i}
\varphi\,d\sigma_i - \frac{1}{2}\int_{\bd{\Omega}_i}
\varphi|\vec{n}_i - \vec{e}_{n + 1}|^2d\sigma_i\Bigr\} =
\int_{\R^n \times \{0\}} \varphi\,dx,$$ so that \eqref{6.6} yields
$$\lim_{i \to \infty} \int_{\bd{\Omega}_i} \varphi|\vec{n}_i -
\vec{e}_{n + 1}|^2d\sigma_i = 0.$$ Taking $\varphi \geq \chi_{B(1,
0)}$ yields the corresponding bound for the integral on
$\bd{\Omega}_i \cap B(1, 0)$.  But
$$\text{av}_{B(1, 0) \cap \bd{\Omega}_i} |\vec{n}_i -
\vec{e}_{n + 1}|^2d\sigma_i = \text{av}_{B(r_i, Q_i)} |\vec{n} -
\vec{e}_{n + 1}|^2d\sigma,$$ and hence $$\ell \leq \varlimsup_{i
\to \infty} \Bigl(\text{av}_{B(r_i, Q_i)} |\vec{n} -
\vec{n}_{B(r_i, Q_i)}|^2d\sigma\Bigr)^{1/2} \leq 2\varlimsup_{i
\to \infty} \Bigl(\text{av}_{B(r_i, Q_i)} |\vec{n} - \vec{e}_{n +
1}|^2d\sigma\Bigr)^{1/2},$$ a contradiction.  This concludes the
proof.

\providecommand{\bysame}{\leavevmode\hbox
to3em{\hrulefill}\thinspace}

\label{lastpage}


\begin{thebibliography}{10}

\bibitem{AC}
H.~W. Alt and L.~A. Caffarelli, Existence and regularity for a
minimum problem with free boundary, \emph{J. Reine Angew. Math.},
325 (1981), 105--144.

\bibitem{D}
B.~Dahlberg, On estimates for harmonic measure, \emph{Arch. Rat.
Mech. Anal.}, 65 (1977), 272--288.

\bibitem{DJ}
G.~David and D.~Jerison, Lipschitz approximation to hypersurfaces,
harmonic measure and singular integrals, \emph{Indiana Univ. Math.
J.}, 39 (1990), 831--845.

\bibitem{DKT}
G.~David, C.~Kenig, and T.~Toro, Asymptotically optimally doubling
measures and Reifenberg flat sets with vanishing constant,
\emph{CPAM}, 54 (2001), 385--449.

\bibitem{Du}
P.~Duren, \emph{The theory of $H^p$ spaces}, Academic Press, New
York, 1970.

\bibitem{EG}
L.~C. Evans and R.~F. Gariepy, \emph{Measure theory and fine
properties of
  functions}, Studies in Advanced Mathematics, CRC Press, 1992.

\bibitem{GCRdeF}
J.~Garc{\'\i}a-Cuerva and J.~L. Rubio~de Francia, \emph{Weighted
norm inequalities and related topics}, Math. Studies, no. 116,
North Holland, 1985.

\bibitem{J}
D.~Jerison, Regularity of the Poisson kernel and free boundary
problems, \emph{Colloquium Mathematicum}, 60--61 (1990), 547--567.

\bibitem{JK1}
D.~Jerison and C.~Kenig, Boundary behavior of harmonic functions
in nontangentially accessible domains, \emph{Adv. in Math.}, 46
(1982), 80--147.

\bibitem{JK2}
\bysame, The logarithm of the Poisson kernel of a $C^1$ domain has
vanishing mean oscillation, \emph{Trans. Amer. Math. Soc.}, 273
(1982), 781--794.

\bibitem{JN}
F.~John and L.~Nirenberg, On functions of bounded mean
oscillation, \emph{Comm. Pure Appl. Math.}, 14 (1961), 415--426.

\bibitem{KL}
M.~V. Keldysh and M.~A. Lavrentiev, Sur la repr\'esentation
conforme des domaines limit\'es par des courbes rectifiables,
\emph{Ann. Sci. Ecole Norm. Sup.}, 54 (1937), 1--38.

\bibitem{KT5}
C.~Kenig and T.~Toro, On the free boundary regularity theorem of
Alt and Caffarelli, preprint.

\bibitem{KT2}
\bysame, Poisson kernel characterization of Reifenberg flat
chord-arc domains, to appear, \emph{Ann. Sci. Ec. Norm. Sup.}

\bibitem{KT3}
\bysame, Harmonic measure on locally flat domains, \emph{Duke
Math. J.}, 87 (1997), 509--551.

\bibitem{KT1}
\bysame, Free boundary regularity for harmonic measures and
Poisson kernels, \emph{Ann. of Math.}, 150 (1999), 369--454.

\bibitem{KoP}
O.~Kowalski and D.~Preiss, Besicovitch-type properties of measures
and submanifolds, \emph{J. Reine Angew. Math.}, 379 (1987),
115--151.

\bibitem{L}
M.~Lavrentiev, Boundary problems in the theory of univalent
functions, \emph{Math Sb. (N. S.) I}, 43 (1936),
  815--844.

\bibitem{P}
Ch.~Pommerenke, On univalent functions, Bloch functions and VMOA,
\emph{Math. Ann.}, 236 (1978), 199--208.

\bibitem{R}
E.~Reifenberg, Solution of the Plateau problem for $m$-dimensional
surfaces of varying topological type, \emph{Acta Math.}, 104
(1960), 1--92.

\bibitem{Se2}
S.~Semmes, Analysis vs. geometry on a class of rectifiable
hypersurfaces, \emph{Indiana Univ. J.}, 39 (1990), 1005--1035.

\bibitem{Se1}
\bysame, Chord-arc surfaces with small constant I,
  \emph{Adv. in Math.}, 85 (1991), 198--293.

\end{thebibliography}
\end{document}